\documentclass{amsart}
\usepackage{amsmath,amssymb,amscd,amsthm}

\newcommand{\C}[1]{\mathcal{#1}}

\newcommand{\B}[1]{\mathbb{#1}}
\newcommand{\xra}[1]{\xrightarrow{\ #1\ }}

\newcommand{\eotimes}[1]{\underset{#1}{\otimes}}
\newcommand{\fotimes}[1]{\underset{#1}{\Box}}

\newtheorem{thm}{Theorem}[section]
\newtheorem{cor}[thm]{Corollary}
\newtheorem{lem}[thm]{Lemma}
\newtheorem{prop}[thm]{Proposition}

\theoremstyle{definition}
\newtheorem{defn}[thm]{Definition}
\newtheorem{rem}[thm]{Remark}

\title{ Hopf Modules and Noncommutative Differential Geometry}
\author{Atabey Kaygun}
\email{{\tt akaygun@uwo.ca}}
\author{Masoud Khalkhali}
\email{{\tt masoud@uwo.ca}}
\address{Department of Mathematics\\
The University of Western Ontario\\
London, Ontario N6A 5B7\\
Canada}
\begin{document}

\begin{abstract}
We define a new algebra of noncommutative differential forms for any
Hopf algebra with an invertible antipode.  We prove that there is a
one to one correspondence between anti-Yetter-Drinfeld modules, which
serve as coefficients for the Hopf cyclic (co)homology, and modules
which admit a flat connection with respect to our differential
calculus.  Thus we show that these coefficient modules can be regarded
as ``flat bundles'' in the sense of Connes' noncommutative
differential geometry.
\end{abstract}

\maketitle

\section{Introduction}

In this paper we define a new algebra of noncommutative differential
forms over any Hopf algebra $H$ with an invertible antipode. The
resulting differential calculus, denoted here by $\mathcal{K}^*(H)$,
is  intimately related to the class of
anti-Yetter-Drinfeld modules over $H$, herein called AYD modules, that
were introduced in \cite{Khalkhali:SaYDModules}.  More precisely, we
show that there is a one to one correspondence between AYD modules
over $H$ and $H$--modules that admit a flat connection with respect to
our differential calculus $\mathcal{K}^*(H)$. In general terms, this
gives a new interpretation of AYD modules as noncommutative analogues
of local systems over the noncommutative space represented by $H$.

The introduction of AYD modules in \cite{Khalkhali:SaYDModules} was
motivated by the problem of finding the largest class of Hopf modules
that could serve as coefficients for the cyclic (co)homology of Hopf
algebras introduced by Connes and Moscovici
\cite{ConnesMoscovici:HopfCyclicCohomology,
ConnesMoscovici:HopfCyclicCohomologyIa,
ConnesMoscovici:HopfCyclicCohomologyII} and its extensions and
ramifications defined in \cite{Khalkhali:DualCyclicHomology,
Khalkhali:InvariantCyclicHomology, Khalkhali:HopfCyclicHomology,
Kaygun:BialgebraCyclicK}.  The question of finding an appropriate
notion of local system, or coefficients, to twist the cyclic homology
of associative algebras is an interesting open problem.  The case of
Hopf cyclic cohomology on the other hand offers an interesting
exception in that it admits coefficients and furthermore, as we shall
prove, these coefficient modules can be regarded as ``flat bundles''
in the sense of Connes' noncommutative differential geometry
\cite{Connes:NonCommutativeGeometry}.

Here is a plan of this paper.  In Section~\ref{Prelims} we recall the
definitions of anti-Yetter-Drinfeld (AYD) and Yetter-Drinfeld (YD)
modules over a Hopf algebra, and we give a characterization of the
category of AYD modules over group algebras as a functor
category. This will serve as our guiding example.  Motivated by this,
in Section~\ref{FlatConnections} we define two natural differential
calculi $\C{K}^*(H)$ and $\widehat{\C{K}}^*(H)$ associated with an
arbitrary Hopf algebra $H$.  Next, we give a complete characterization
of the category of AYD and YD modules over a Hopf algebra $H$ as
modules admitting flat connections over these differential calculi
$\C{K}^*(H)$ and $\widehat{\C{K}}^*(H)$, respectively.  In the same
section, we also investigate tensor products of modules admitting flat
connections.  In Section~\ref{EquivariantCalculi}, we place the
differential calculi $\C{K}^*(H)$ and $\widehat{\C{K}}^*(H)$ in a much
larger class of differential calculi, unifying the results proved
separately before. This last section is motivated by the recent paper
\cite{PanaiteStaic:GeneralizedAYD}.  After this paper was posted, T.
Brzezinski informed us that he can also derive our differential
calculus and flatness results from the theory of corings
\cite{BrzezinskiWisbauer:ComodulesCorings}. This would involve
realizing AYD modules via a specific entwining structure, and using
the general machinery of \cite{BrzezinskiWisbauer:ComodulesCorings}
that cast such structures as flat modules.

In this paper we work over a field $k$ of an arbitrary characteristic.
All of our results, however, are valid over an arbitrary unital
commutative ground ring $k$ if the objects involved (coalgebras,
bialgebras, Hopf algebras and their (co)modules) are flat
$k$--modules.  We will assume $H$ is a Hopf algebra with an invertible
antipode and all of our $H$--modules and comodules are left modules
and comodules unless explicitly indicated otherwise.

Acknowledgement: M.K. would like to thank professor Alain Connes for
suggesting the problem of interpreting AYD modules as flat modules in
the sense of noncommutative differential geometry and for enlightening
discussions at an early stage of this work.

\section{Preliminaries and a guiding example}\label{Prelims}

Recall from \cite{Khalkhali:SaYDModules} that a $k$--module $X$ is
called a left-left anti-Yetter-Drinfeld (AYD, for short) module if (i)
$X$ is a left $H$--module, (ii) $X$ is a left $H$--comodule, and (iii)
one has a compatibility condition between the $H$--module and comodule
structure on $X$ in the sense that:
\begin{align}\label{AYD}
  (hx)_{(-1)}\otimes(hx)_{(0)}
   = h_{(1)}x_{(-1)}S^{-1}(h_{(3)})\otimes h_{(2)}x_{(0)}
\end{align}
for any $h\in H$ and $x\in X$.  The AYD condition~(\ref{AYD}) should
be compared with the Yetter-Drinfeld (YD) compatibility condition:
\begin{align}\label{YD}
  (hx)_{(-1)}\otimes (hx)_{(0)}
   = h_{(1)}x_{(-1)}S(h_{(3)})\otimes h_{(2)}x_{(0)}
\end{align}
for any $h\in H$ and $x\in X$.

A morphism of AYD modules $X\to Y$ is simply a $k$-linear map $X
\to Y$ compatible with the $H$--action and coaction. The resulting
category of AYD modules over $H$ is an abelian category. There is a
similar statement for YD modules.

To pass from (\ref{AYD}) to (\ref{YD}) one replaces $S^{-1}$ by $S$.
This is, however, a nontrivial operation since the categories of AYD
and YD modules over $H$ are very different in general. For example,
the former is not a monoidal category in a natural way
\cite{Khalkhali:SaYDModules} while the latter is always monoidal. If
$S^2=id_H$, in particular when $H$ is commutative or cocommutative,
then these categories obviously coincide.

To understand these compatibility conditions better we proceed as
follows.  Let $X$ be a left $H$ module. We define a left $H$--action on
$H\otimes X$ by letting
\begin{align}\label{AYDact}
h (g\otimes x): =h_{(1)}gS^{-1}(h_{(3)})\otimes h_{(2)}x
\end{align}
for any $h\in H$ and $g\otimes x \in H\otimes X$.

\begin{lem}
The formula given in (\ref{AYDact}) defines a left $H$--module
structure on $H\otimes X$.  Moreover, an $H$--module/comodule $X$ is
an anti-Yetter-Drinfeld module iff its comodule structure map
$\rho_X:X\to H\otimes X$ is a morphism of $H$--modules.
\end{lem}

There is of course a similar characterization of YD modules. The left
action (\ref{AYDact}) should simply be replaced by the left action
\begin{align}\label{YDact}
h (g\otimes x):=h_{(1)}gS(h_{(3)})\otimes h_{(2)}x
\end{align}

Let us give a characterization of AYD modules in a concrete example.
Let $G$ be a, not necessarily finite, discrete group and let $H=k[G]$
be its groups algebra over $k$ with its standard Hopf algebra
structure, i.e. $\Delta(g)=g\otimes g$ and $S(g)=g^{-1}$ for all $g\in
G$.  We define a groupoid $G\ltimes G$ whose set of morphisms is
$G\times G$ and its set of objects is $G$.  Its source and target maps
$s: G\ltimes G\to G$ and $t:G\ltimes G\to G$ are defined by
\begin{align*}
  s(g,g') &= g  & t(g,g') &= gg'g^{-1}
\end{align*}
for any $(g,g')$ in $G\ltimes G$. It is easily seen that under group
multiplication as its composition, $G\ltimes G$ is a groupoid.

\begin{prop}\label{BundlesOverG}
  The category of Yetter-Drinfeld modules over $k[G]$ is
  isomorphic to the category of functors from $G\ltimes G$ into the
  category of $k$--modules.
\end{prop}

\begin{proof}
  Let $M$ be a  $k[G]$--module/comodule.  Denote its  structure morphisms
  by $\mu:k[G] \otimes M\to M$ and $\rho:M\to k[G] \otimes M$.
  Since we assumed $k$ is a field, $M$ has a basis of the form
  $\{e^i\}_{i\in I}$ for some index set $I$.  Since $M$ is a $k[G]
  $--comodule one has
  \begin{align*}
    e^i_{(-1)}\otimes e^i_{(0)} = \sum_{j\in I}\sum_{g\in G}
    c^i_{j,g}(g\otimes e^j)
  \end{align*}
  where only finitely many $c_{j,g}$ is non-zero.  One can chose a
  basis $\{m^\lambda\}_{\lambda\in\Lambda}$ for $M$ such that
  \begin{align*}
    m^\lambda_{(-1)}\otimes m^\lambda_{(0)}
    = \sum_\alpha c_{\lambda,\alpha} (g_\lambda\otimes m^\alpha)
  \end{align*}
  and since all comodules are counital and $k[G]$ has a counit
  $\varepsilon(g)=1$ for any $g\in G$ we see that
  \begin{align*}
    m^\lambda = \sum_\alpha c_{\lambda,\alpha} m^\alpha
  \end{align*}
  implying $c_{\lambda,\alpha}$ is uniformly zero except
  $c_{\lambda,\lambda}$ which is 1. In other words, one can split $M$
  as $\bigoplus_{g\in G} M_g$ such that $\rho(x)= g\otimes x$ for any
  $x\in M_g$.  Now assume $M$ is an AYD module.  Then since
  \begin{align*}
    (hx)_{(-1)}\otimes (hx)_{(0)} = hgh^{-1}\otimes hx
  \end{align*}
  for any $x\in M_g$ and $h\in G$ one can see that
  $L_h:M_g\to M_{hgh^{-1}}$ where $L_h$ is the $k$--vector space
  endomorphism of $M$ coming from the left action of $h$.  This
  observation implies that the category of AYD modules over $k[G]$
  and the category of functors from $G\ltimes G$ into the category of
  $k$--modules are isomorphic.
\end{proof}

\section{Differential calculi and flat connections}\label{FlatConnections}

Our next goal is to find a noncommutative analogue of
Proposition~\ref{BundlesOverG}.  To this end, we will replace the
groupoid $G\ltimes G$ by a differential calculus $\C{K}^*(H)$
naturally defined for any Hopf algebra $H$.  The right analogue of
representations of the groupoid $G\ltimes G$ will be $H$--modules
admitting flat connections with respect to the differential calculus
$\C{K}^*(H)$.

Let us first recall basic notions of connection and curvature in the
noncommutative setting from \cite{Connes:NonCommutativeGeometry}.  Let
$A$ be a $k$--algebra.  A differential calculus over $A$ is a
differential graded $k$--algebra $(\Omega^*,d)$ endowed with a
morphism of algebras $\rho:A\to \Omega^0$.  The differential $d$ is
assumed to have degree one.  Since in our main examples we have
$\Omega^0=A$ and $\rho=id$, in the following we assume this is the
case.

Assume $M$ is a left $A$--module.  A morphism of $k$--modules
$\nabla:M\to \Omega^1\eotimes{A}M$ is called a connection with
respect to the differential calculus $(\Omega^*,d)$ if one has a
Leibniz rule of the form
\begin{align*}
  \nabla(am) = a\nabla(m) +  d(a)\eotimes{A}m
\end{align*}
for any $m\in M$ and $a\in A$.  Given any connection $\nabla$ on $M$,
there is a unique extension of $\nabla$ to a map $\widehat{\nabla}:
\Omega^*\eotimes{A}M\to\Omega^*\eotimes{A}M$ satisfying a graded
Leibniz rule.  It is given by
\begin{align*}
  \widehat{\nabla}(\omega\otimes m)
   = d(\omega)\eotimes{A}m + (-1)^{|\omega|}\omega\nabla(m)
\end{align*}
for any $m\in M$ and $\omega\in\Omega^*$.  A connection
$\nabla:M\to \Omega^1\eotimes{A}M$ is called flat if its curvature
$R: =\widehat{\nabla}^2=0$.  This is equivalent to saying
$\Omega^*\eotimes{A}M$ is a differential graded $\Omega^*$--module
with the extended differential $\widehat{\nabla}$.

The following general definition will be useful in the rest of this
paper.

\begin{defn}
  Let $X_0,\ldots,X_n$ be a finite set of $H$--bimodules.  We define
  an $H$--bimodule structure on the $k$--module
  $X_0\otimes\cdots\otimes X_n$ by
  \begin{align*}
    h(x^0\otimes\cdots\otimes x^n)
     = & h_{(1)}x^0 S^{-1}(h_{(2n+1)})\otimes\cdots\otimes
         h_{(n)} x^{n-1} S^{-1}(h_{(n+2)})\otimes h_{(n+1)} x^n,\\
    (x^0\otimes\cdots\otimes x^n)h
     = & x^0\otimes\cdots\otimes x^{n-1}\otimes x^nh.
  \end{align*}
  for any $h\in H$ and $(x^0\otimes\cdots\otimes x^n)\in
  X_0\otimes\cdots\otimes X_n$.  Checking the bimodule conditions is
  straightforward. We denote this bimodule by $X_0\oslash\cdots\oslash
  X_n$.
\end{defn}

\begin{rem}
  We should  remark that $X\oslash Y$   is not a monoidal product.
  In other words given any three $H$--bimodules $X$, $Y$, and $Z$, then
   $(X\oslash Y)\oslash Z$ and $X\oslash (Y\oslash Z)$ are
  not isomorphic as left $H$--modules unless $H$ is cocommutative.
\end{rem}

\begin{defn}
  For each $n\geq 0$, let $\C{K}^n(H)= H^{\oslash n+1}$.
  We define a differential $d: \C{K}^n(H) \to
  \C{K}^{n+1}(H)$ by
  \begin{align*}
    d(h^0\otimes\cdots\otimes h^n)
     = &  - (1\otimes h^0\otimes\cdots\otimes h^n)
          + \sum_{j=0}^{n-1}(-1)^j(h^0\cdots\otimes
                h^j_{(1)}\otimes h^j_{(2)}\otimes\cdots\otimes h^n)\\
       & + (-1)^n(h^0\otimes\cdots\otimes h^{n-1}\otimes
                  h^n_{(1)}S^{-1}(h^n_{(3)})\otimes h^n_{(2)}).
  \end{align*}
    We also
  define an associative graded product structure by
  \begin{align*}
    (x^0\otimes\cdots\otimes x^n) & (y^0\otimes\cdots\otimes y^m)\\
    = &  x^0\otimes\cdots\otimes x^{n-1}\otimes
         x^n_{(1)}y^0S^{-1}(h_{(2m+1)})\otimes\cdots\otimes
         x^n_{(m)}y^{m-1}S^{-1}(x^n_{(m+2)})\otimes x^n_{(m+1)}y^m
  \end{align*}
  for any $(x^0\otimes\cdots\otimes x^n)$ in $\C{K}^n(H)$ and
  $(y^0\otimes\cdots\otimes y^m)$ in $\C{K}^m(H)$.
\end{defn}

\begin{prop}\label{AYDCalculus}
  $\C{K}^*(H)$ is a differential graded $k$--algebra.
\end{prop}

\begin{proof}
  For any $x\in\C{K}^0(H)$ one has
  \begin{align*}
    d(x) = - (1\otimes x) + (x_{(1)}S^{-1}(x_{(3)})\otimes x_{(2)})
         = [x,(1\otimes 1)],
  \end{align*}
  and for $(y\otimes 1)$ in $\C{K}^1(H)$ and $x\in\C{K}^0(H)$ we see
  \begin{align*}
    d(x(y\otimes 1))
    = & d(x_{(1)}yS^{-1}(x_{(3)})\otimes x_{(2)})\\
    = & - (1\otimes x_{(1)}yS^{-1}(x_{(3)})\otimes x_{(2)})
        + (x_{(1)}y_{(1)}S^{-1}(x_{(5)})\otimes x_{(2)}y_{(2)}S^{-1}(x_{(4)})\otimes x_{(3)})\\
      & - (x_{(1)}yS^{-1}(x_{(5)})\otimes x_{(2)}S^{-1}(x_{(4)})\otimes x_{(3)})\\
    %= & - (1\otimes x_{(1)}yS^{-1}(x_{(3)})\otimes x_{(2)})
    %    + (x_{(1)}S^{-1}(x_{(5)})\otimes x_{(2)}y S^{-1}(x_{(4)})\otimes x_{(3)})\\
    %  & - (x_{(1)}S^{-1}(x_{(5)})\otimes x_{(2)}y S^{-1}(x_{(4)})\otimes x_{(3)})
    %    + (x_{(1)}y_{(1)}S^{-1}(x_{(5)})\otimes x_{(2)}y_{(2)}S^{-1}(x_{(4)})\otimes x_{(3)})\\
    %  & -(x_{(1)}yS^{-1}(x_{(5)})\otimes x_{(2)}S^{-1}(x_{(4)})\otimes x_{(3)})\\
    = & d(x)(y\otimes 1) + x d(y\otimes 1)
  \end{align*}
  We also see for $(x\otimes y)$ in $\C{K}^1(H)$ the we have
  \begin{align*}
    d((x\otimes 1)y) & = d(x\otimes y)
    =   - (1\otimes x\otimes y)
        + (x_{(1)}\otimes x_{(2)}\otimes y)
        - (x\otimes y_{(1)}S^{-1}(y_{(3)})\otimes y_{(2)})\\
    = & - (1\otimes x\otimes 1)y
        + (x_{(1)}\otimes x_{(2)}\otimes 1) y
        - (x\otimes 1\otimes 1)y
        + (x\otimes 1)(1\otimes y)\\
      & - (x\otimes 1)(y_{(1)}S^{-1}(y_{(3)})\otimes y_{(2)})\\
    = & d(x\otimes 1)y - (x\otimes 1)d(y).
  \end{align*}
  Note that with the product structure on $\C{K}^*(H)$ one
  has
  \begin{align*}
    (x^0\otimes\cdots\otimes x^n)
    = (x^0\otimes 1)\cdots(x^{n-2}\otimes 1)(x^{n-1}\otimes 1)x^n
  \end{align*}
  for any $x^0\otimes\cdots\otimes x^n$ in  $\C{K}^n(H)$.  Now,
  one can inductively show that
  \begin{align*}
    d(\Psi\Phi)
    = & d(\Psi)\Phi + (-1)^{|\Psi|}\Psi d(\Phi)
  \end{align*}
  for any $\Psi$ and $\Phi$ in $\C{K}^*(H)$.  Since the algebra is
  generated by degree zero and degree one terms, all that remains is
  to show that for all $x\in H$ we have $d^2(x)=0$ and $d^2(x\otimes
  1)=0$.  For the first assertion we see that
  \begin{align*}
    d^2(x) & = - d(1\otimes x) + d(x_{(1)}S^{-1}(x_{(3)})\otimes x_{(2)})\\
     = & (1\otimes 1\otimes x) - (1\otimes 1\otimes x)
        + (1\otimes x_{(1)}S^{-1}(x_{(3)})\otimes x_{(2)})
        - (1\otimes x_{(1)}S^{-1}(x_{(3)})\otimes x_{(2)})\\
       &+ (x_{(1)(1)}S^{-1}(x_{(3)(2)})\otimes x_{(1)(2)}S^{-1}(x_{(3)(1)})\otimes x_{(2)})
        - (x_{(1)}S^{-1}(x_{(3)})\otimes x_{(2)(1)}S^{-1}(x_{(2)(3)})\otimes x_{(2)(2)})\\
     = & 0
  \end{align*}
  for any $x\in H$.  For the second assertion we see
  \begin{align*}
    d^2& (x\otimes 1)
     =  - d(1\otimes x\otimes 1)
         + d(x_{(1)}\otimes x_{(2)}\otimes 1)
         - d(x\otimes 1\otimes 1)\\
     = & - (1\otimes 1\otimes x\otimes 1)
         + (1\otimes 1\otimes x\otimes 1)
         - (1\otimes x_{(1)}\otimes x_{(2)}\otimes 1)
         + (1\otimes x\otimes 1\otimes 1)\\
       & - (1\otimes x_{(1)}\otimes x_{(2)}\otimes 1)
         + (x_{(1)}\otimes x_{(2)}\otimes x_{(3)}\otimes 1)
         - (x_{(1)}\otimes x_{(2)}\otimes x_{(3)}\otimes 1)
         + (x_{(1)}\otimes x_{(2)}\otimes 1\otimes 1)\\
       & + (1\otimes x\otimes 1\otimes 1)
         - (x_{(1)}\otimes x_{(2)}\otimes 1\otimes 1)
         + (x\otimes 1\otimes 1\otimes 1)
         - (x\otimes 1\otimes 1\otimes 1)\\
     = & 0
  \end{align*}
  for any $(x\otimes 1)$ in $\C{K}^*(H)$.  The result follows.
\end{proof}

Note that the calculus $\C{K}^*(H)$ is determined by (i) the
$H$--bimodule $\C{K}^1(H) = H\otimes H$ (ii) the differential
$d_0:H\to H\otimes H$ and $d_1:H\otimes H\to H\otimes H\otimes H$
and (iii) the Leibniz rule $d(\Psi\Phi) = d(\Psi)\Phi +
(-1)^{|\Psi|}\Psi d(\Phi)$.

Recall  that for any coassociative coalgebra $C$, the
category of $C$--bicomodules has a monoidal product called cotensor
product, which is denoted by $\fotimes{C}$.  The cotensor product is
left exact and its right derived functors are denoted by ${\rm
Cotor}_C^*(\cdot,\cdot)$.

Now we can identify the homology of the calculus $\C{K}^*(H)$ as
follows:
\begin{prop}
  $H_*(\C{K}^*(H))$ is isomorphic to ${\rm
  Cotor}_H^*(k,H^{coad})$ where $k$ is considered as an $H$--comodule
  via the unit and $H^{coad}$ is the coadjoint corepresentation over
  $H$, i.e. $\rho^{coad}(h)=h_{(1)}S^{-1}(h_{(3)})\otimes h_{(2)}$ for
  any $h\in H$.
\end{prop}

\begin{thm}\label{AYDModules}
  The category of AYD modules over $H$ is isomorphic to the category
  of $H$--modules admitting a flat connection with respect to  the
  differential calculus $\C{K}^*(H)$.
\end{thm}

\begin{proof}
  Assume $M$ is a $H$--module which admits a morphism of $k$--modules
  of the form $\nabla:M\to \C{K}^1(H)\eotimes{H}M\cong H\otimes M$.
  Define $\rho_M(m)=\nabla(m)+(1\otimes m)$ and denote $\rho_M(m)$ by
  $(m_{(-1)}\otimes m_{(0)})$ for any $m\in M$.  First we see that
  \begin{align*}
    \nabla(hm) = & (hm)_{(-1)}\otimes (hm)_{(0)} - (1\otimes hm)
  \end{align*}
  and also
  \begin{align*}
  d(h)\eotimes{H}m + h\nabla(m)
      = & - (1\otimes hm)
          + (h_{(1)}S^{-1}(h_{(3)})\otimes h_{(2)}m)\\
        & + (h_{(1)}m_{(-1)}S^{-1}(h_{(3)})\otimes h_{(2)}m_{(0)})
        - (h_{(1)}S^{-1}(h_{(3)})\otimes h_{(2)}m)\\
      = & (h_{(1)}m_{(-1)}S^{-1}(h_{(3)})\otimes h_{(2)}m_{(0)})
          - (1\otimes hm)
  \end{align*}
  for any $h\in H$ and $m\in M$.  This means $\nabla$ is a connection
  iff the $H$--module $M$ together with $\rho_X:M\to H\otimes M$
   satisfy the AYD condition.  The flatness condition will hold iff
  for any $m\in M$ one has
  \begin{align*}
    \widehat{\nabla}^2(m)
      = & d(m_{(-1)}\otimes 1)\eotimes{H} m_{(0)}
         - (m_{(-1)}\otimes 1)\nabla(m_{(0)})
         - d(1\otimes 1)m + (1\otimes 1)\nabla(m)\\
      %= & - (1\otimes m_{(-1)}\otimes m_{(0)})
      %    + (m_{(-1)(1)}\otimes m_{(-1)(2)}\otimes m_{(0)})
      %    - (m_{(-1)}\otimes 1\otimes m_{(0)})\\
      %  & - (m_{(-1)}\otimes m_{(0)(-1)}\otimes m_{(0)(0)})
      %    + (m_{(-1)}\otimes 1\otimes m_{(0)})\\
      %  & + (1\otimes 1\otimes m) + (1\otimes m_{(-1)}\otimes m_{(0)})
      %    - (1\otimes 1\otimes m)\\
      = & (m_{(-1)(1)}\otimes m_{(-1)(2)}\otimes m_{(0)})
          - (m_{(-1)}\otimes m_{(0)(-1)}\otimes m_{(0)(0)})
      = 0,
  \end{align*}
  meaning $\nabla$ is flat iff $\rho_M:M\to H\otimes M$ defines a
  coassociative coaction of $H$ on $M$.
\end{proof}

\begin{defn}
  Let $X$ be an AYD module over $H$.  Define $\C{K}^*(H,X)$
  as the graded $k$--module $\C{K}^*(H)\eotimes{H}X$ equipped with the
  connection $\nabla_X(x) = \rho_X(x) - (1\otimes x)$ as its
  differential.
\end{defn}

Recall from \cite{Khalkhali:SaYDModules} that an $H$--module/comodule
$X$ is called stable if the composition $X\xra{\rho_X}H\otimes
X\xra{\mu_X}X$ is $id_X$ where $\rho_X$ and $\mu_X$ denote the
$H$--comodule and $H$--module structure maps respectively.  Explicitly
one has $x_{(-1)}x_{(0)}=x$ for any $x\in X$.

\begin{thm}
  For an arbitrary AYD module $X$ one has $H_*(\C{K}^*(H,X))\cong {\rm
  Cotor}_H^*(k,X)$.  Moreover, if $X$ is also stable, then
  $\C{K}^*(H,X)$ is isomorphic (as differential graded $k$--modules)
  to the Hochschild complex of the Hopf-cyclic complex of the Hopf
  algebra $H$ with coefficients in $X$.
\end{thm}

\begin{proof}
  The first part of the Theorem follows from the observation that
  $\C{K}^*(H,X)$, viewed just as a differential graded $k$--module, is
  really $\C{B}_*(k,H,X)$ the two sided cobar complex of the coalgebra
  $H$ with coefficients in $H$--comodules $k$ and $X$.  The second
  assertion follows from Remark~3.13 and Theorem~3.14 of
  \cite{Kaygun:BialgebraCyclicK}.
\end{proof}

\begin{prop}
  Let $X$ be an AYD module over $H$.  If $H$ is cocommutative, then
  $\C{K}^*(H,X)$ is a differential graded left $H$--module with
  respect to the AYD module structure on $\C{K}^*(H,X)$.
\end{prop}

Instead of the AYD condition, one can consider the YD condition and form a
differential calculus $\widehat{\C{K}}^*(H)$ using the YD condition.

\begin{defn}
  As before, assume $H$ is a Hopf algebra, but this time we do not
  require the antipode to be invertible.  We define a new differential
  calculus $\widehat{\C{K}}^*(H)$ over $H$ as follows: let
  $\widehat{\C{K}}^n(H)=H^{\otimes n+1}$ and define the differentials
  as
  \begin{align*}
    d(x^0\otimes\cdots\otimes x^n)
    = & -(1\otimes x^0\otimes\cdots\otimes x^n)
      + \sum_{j=0}^{n-1}(-1)^j(x^0\otimes \cdots\otimes x^j_{(1)}\otimes x^j_{(2)}\otimes\cdots x^n)\\
      & + (-1)^n(x^0\otimes\cdots\otimes x^{n-1}\otimes x^n_{(1)}S(x^n_{(3)})\otimes x^n_{(2)})
  \end{align*}
  for any $x^0\otimes\cdots\otimes x^n$ in $\widehat{\C{K}}^n(H)$.
  The multiplication is defined as
  \begin{align*}
    (x^0\otimes\cdots\otimes x^n) & (y^0\otimes\cdots\otimes y^m)\\
    = & x^0\otimes\cdots\otimes x^{n-1}\otimes x^n_{(1)}y^0S(x_{(2m+1)})
         \otimes\cdots\otimes x^n_{(m)}y^{m-1}S(x^n_{(m+1)})\otimes x^n_{(m+1)}y^m
  \end{align*}
  for any $x^0\otimes\cdots\otimes x^n$ and $y^0\otimes\cdots\otimes
  y^m$ in $\widehat{\C{K}}^*(H)$.
\end{defn}

The proofs of the following facts are similar to  the corresponding
statements for the differential calculus $\C{K}^*(H)$ and AYD modules.

\begin{prop}
  $\widehat{\C{K}}^*(H)$ is a differential graded $k$--algebra.
\end{prop}

\begin{thm}\label{YDModules}
  The category of YD modules over $H$ is isomorphic to the category of
  $H$--modules admitting a flat connection with respect to the
  differential
  calculus $\widehat{\C{K}}^*(H)$.
\end{thm}

\begin{defn}
  Let $X$ be a Yetter-Drinfeld module $X$ over $H$ with the structure
  morphisms $\mu_X:H\otimes X\to X$ and $\rho_X:X\to H\otimes X$.
  Define $\widehat{\C{K}}^*(H,X)$ as the (differential) graded
  $k$--module $\widehat{\C{K}}^*(H)\eotimes{H}X$ with the connection
  $\nabla_X(x) = \rho_X(x) - (1\otimes x)$ defined for any $x\in X$.
\end{defn}

\begin{prop}
  For an arbitrary YD module $X$ one has
  $H_*(\widehat{\C{K}}^*(H,X))\cong {\rm Cotor}_H^*(k,X)$.
\end{prop}

\begin{prop}
  Assume $X$ is an arbitrary YD module over $H$.  If $H$ is
  cocommutative, then $\widehat{\C{K}}^*(H,X)$ is a differential
  graded left $H$--module with respect to the YD module structure on
  $\widehat{\C{K}}^*(H)$.
\end{prop}

Our next goal is to study tensor products of AYD and YD modules in our
noncommutative differential geometric setup. It is well known that
the tensor product of two flat vector bundles over a manifold is again
a flat bundle.  Moreover, from the resulting monoidal, in fact
Tannakian, category one can recover the fundamental group of the base
manifold.  The situation in the noncommutative case is of course far
more complicated and we only have some vestiges of this theory.

Assume $H$ is a Hopf algebra with an invertible antipode.  Let $X$ be
a $H$--module admitting a flat connection $\nabla$ with respect to  the calculus
$\C{K}^*(H)$.  We define a switch morphism $\sigma: X\otimes
H\to H\otimes X$ by letting
\begin{align*}
  \sigma(x\otimes h) = x_{\{-1\}} h\otimes x_{\{0\}}+h\otimes x,
\end{align*}
for any $x\otimes h$ in $X\otimes H$ where we used a Sweedler notation
to denote the connection: $\nabla(x)=x_{\{-1\}} \otimes
x_{\{0\}}$. Note that $\sigma$ is a perturbation of the standard
switch map.

\begin{prop}
  Let $X$ and $X'$ be two $H$--modules with  flat connections
  $\nabla_X$ and $\nabla_{X'}$ with respect to the differential
  calculi $\widehat{\C{K}}^*(H)$ and $\C{K}^*(H)$ respectively.  Then
  $\nabla_{X\otimes X'}:X\otimes X'\to H\otimes X\otimes X'$ given
  by
  \begin{align*}
    \nabla_{X\otimes X'}(x\otimes x')=\nabla_X(x)\otimes x'
      + (\sigma\otimes id_{X'})\left(x\otimes \nabla_{X'}(x')\right)
  \end{align*}
  for any $x\otimes x'$ in $X\otimes X'$ defines a flat connection
  on the $H$--module $X\otimes X'$ with respect to $\C{K}^*(H)$.
\end{prop}

\begin{proof}
Recall from Theorem~\ref{AYDModules} and Theorem~\ref{YDModules} that
the category of AYD (resp. YD) modules and the category of $H$--modules
admitting flat connections with respect to the differential calculus
$\C{K}^*(H)$ (resp. $\widehat{\C{K}}^*(H)$) are isomorphic.  Then
given two $H$--modules with flat connections $(X,\nabla_X)$ over
$\widehat{\C{K}}^*(H)$ and $(X',\nabla_{X'})$ over $\C{K}^*(H)$ one
can extract $H$--comodule structures by letting $x_{(-1)}\otimes
x_{(0)}:=\rho_X(x):=\nabla_X(x)+(1\otimes x)$ and $x'_{(-1)}\otimes
x'_{(0)}:=\rho_{X'}(x'):=\nabla_{X'}(x')+(1\otimes x')$.  Then we get
\begin{align*}
\nabla_{X\otimes X'}(x\otimes x')
  = & \nabla_X(x)\otimes x'
      + (\sigma\otimes id_{X'})\left(x\otimes \nabla_{X'}(x')\right)\\
  = & - (1\otimes x\otimes x') + (x_{(-1)}\otimes x_{(0)}\otimes x')
      - (x_{(-1)}\otimes x_{(0)}\otimes x')
      + (x_{(-1)} x'_{(-1)}\otimes x_{(0)}\otimes x'_{(0)})\\
  = & - (1\otimes x\otimes x')
      + (x_{(-1)} x'_{(-1)}\otimes x_{(0)}\otimes x'_{(0)})
\end{align*}
for any $(x\otimes x')$ in $X\otimes X'$.  One can easily check that
$\rho_{X\otimes X'}(x\otimes x')=x_{(-1)}x'_{(-1)}\otimes
x_{(0)}\otimes x'_{(0)}$ is a coassociative $H$--coaction on $X\otimes
X'$.  Moreover, for any $h\in H$ and $x\otimes x'$ in $X\otimes X'$
we have
\begin{align*}
  \rho_{X\otimes X'}(h(x\otimes x'))
  = & \rho_{X\otimes X'}(h_{(1)}(x)\otimes h_{(2)}(x'))\\
  = & h_{(1)}x_{(-1)}S(h_{(3)})h_{(4)}x'_{(-1)}S^{-1}(h_{(6)})\otimes
        h_{(2)}(x_{(0)})\otimes h_{(5)}x'_{(0)}\\
  = & h_{(1)}x_{(-1)}x'_{(-1)}S^{-1}(h_{(3)})\otimes
        h_{(2)}(x_{(0)}\otimes x'_{(0)})
\end{align*}
meaning $X\otimes X'$ is an AYD module over $H$.  In other words
$\nabla_{X\otimes X'}$ is a flat connection on $X\otimes X'$ with
respect to the differential calculus $\C{K}^*(H)$.  The result
follows.
\end{proof}

\section{Equivariant differential calculi}\label{EquivariantCalculi}

The concepts of AYD and YD modules have recently been extended in
\cite{PanaiteStaic:GeneralizedAYD}. In this section we give a further
extension of this new class of modules and show that they can be
interpreted as modules admitting flat connections with respect to a
differential calculus.  In this section we assume $B$ is a bialgebra
and $B^{op, cop}$ is the bialgebra $B$ with the opposite
multiplication and comultiplication.  Assume $\alpha:B\to B$ and
$\beta:B\to B^{op,cop}$ are two morphisms bialgebras.  Also in this
section we fix a $B$--bimodule coalgebra $C$.  In other words $C$ is a
$B$--bimodule and the comultiplication $\Delta_C:C\to C\otimes C$ is a
morphism of $B$--bimodules where we think of $C\otimes C$ as a
$B$--bimodule via the diagonal action of $B$.  Equivalently, one has
\begin{align*}
  (bcb')_{(1)} = b_{(1)}c_{(1)}b'_{(1)}\otimes b_{(2)}c_{(2)}b'_{(2)}
\end{align*}
for any $b,b'\in B$ and $c\in C$.  We also assume $C$ has a grouplike
element via a coalgebra morphism $\B{I}:k\to C$.  We do not impose
any condition on $\varepsilon(\B{I})$.

\begin{defn}
  Define a graded $B$--bimodule $\C{K}_{(\alpha,\beta)}^*(C,B)$ by
  letting $\C{K}_{(\alpha,\beta)}^n(C,B) = C^{\otimes n}\otimes B$.
  The right action is defined by the right regular representation of
  $B$ on itself, i.e.
  \begin{align*}
    (c^1\otimes\cdots\otimes c^n\otimes b')b
    = c^1\otimes\cdots\otimes c^n\otimes b' b
  \end{align*}
  and the left action is defined by
  \begin{align*}
    b(c^1\otimes\cdots\otimes c^n\otimes b')
    = \alpha(b_{(1)}) c^1 \beta(b_{(2n+1)})\otimes\cdots\otimes
      \alpha(b_{(n)}) c^n \beta(b_{(n+2)})
      \otimes b_{(n+1)} b'
  \end{align*}
  for any $c^1\otimes\cdots\otimes c^n\otimes b'$ in
  $\C{K}_{(\alpha,\beta)}^n(C,B)$ and $b\in B$.
\end{defn}

\begin{prop}
  $\C{K}_{(\alpha,\beta)}^*(C,B)$ is a differential graded
  $k$--algebra.
\end{prop}

\begin{proof}
  First we define a product structure.  We let
  \begin{align*}
    (c^1\otimes\cdots\otimes & c^n\otimes b)
    (c^{n+1}\otimes\cdots\otimes c^{n+m}\otimes b')\\
    = &  c^1\otimes\cdots\otimes c^n\otimes
        \alpha(b_{(1)}) c^{n+1} \beta(b_{(2m+1)})\otimes\cdots\otimes
        \alpha(b_{(m)}) c^{n+m} \beta(b_{(m+2)})\otimes b_{(m+1)}b'
  \end{align*}
  for any $c^1\otimes\cdots\otimes c^n\otimes b$ and
  $c^{n+1}\otimes\cdots\otimes c^{n+m}\otimes b'$ in
  $\C{K}_{(\alpha,\beta)}^*(C,B)$.  Note that
  \begin{align*}
    c^1\otimes\cdots\otimes c^n\otimes b
     = (c^1\otimes 1)\cdots(c^n\otimes 1)b.
  \end{align*}
  So it is enough to check associativity only for degree 1 terms.
  Then
  \begin{align*}
    ((c^1\otimes b^1)(c^2\otimes b^2))(c^3\otimes b^3)
     = & (c^1\otimes\alpha(b^1_{(1)}) c^2 \beta(b^1_{(3)})\otimes b^1_{(2)} b^2)
         (c^3\otimes b^3)\\
%     = & c^1\otimes \alpha(b^1_{(1)}) c^2 \beta(b^1_{(3)})\otimes
%         \alpha(b^1_{(2)(1)})\alpha(b^2_{(1)}) c^3 \beta(b^2_{(3)})\beta(b_{(2)(3)})\otimes
%         b^1_{(2)(2)}b^2_{(2)} b^3\\
     = & c^1\otimes \alpha(b^1_{(1)}) c^2 \beta(b^1_{(5)})\otimes
         \alpha(b^1_{(2)})\alpha(b^2_{(1)}) c^3 \beta(b^2_{(3)})\beta(b^1_{(4)})\otimes
         b^1_{(3)}b^2_{(2)} b^3\\
     = & (c^1\otimes b^1)(c^2\otimes \alpha(b^2_{(1)}) c^3 \beta(b^2_{(3)})\otimes
          b^2_{(2)} b^3)\\
     = & (c^1\otimes b^1)((c^2\otimes b^2)(c^3\otimes b^3))
  \end{align*}
  for any $(c^i\otimes b^i)$ for $i=1,2,3$ as we wanted to show.
  Define the differentials as
  \begin{align*}
    d(c^1\otimes\cdots\otimes c^n\otimes b)
    = & -(\B{I}\otimes c^1\otimes\cdots\otimes c^n\otimes b)
        + \sum_{j=1}^n (-1)^{j-1}(c^1\otimes\cdots\otimes c^j_{(1)}\otimes c^j_{(2)}
         \otimes\cdots\otimes c^n\otimes b)\\
      & (-1)^n (c^1\otimes\cdots\otimes c^n\otimes
         \alpha(b_{(1)})\B{I}\beta(b_{(3)})\otimes b_{(2)})
  \end{align*}
  for any $(c^1\otimes\cdots\otimes c^n\otimes b)$ in
  $\C{K}_{(\alpha,\beta)}^*(C,B)$ and one can check that
  \begin{align*}
    d(b) & = -(\B{I}\otimes b) + (\alpha(b_{(1)})\B{I}\beta(b_{(3)})\otimes b_{(2)}) &
    d(c\otimes 1) = & -(\B{I}\otimes c\otimes 1) + (c_{(1)}\otimes c_{(2)}\otimes 1)
                      -(c\otimes\B{I}\otimes 1)
  \end{align*}
  for any $b\in B$, $c\in C$.  In order to prove that
  $\C{K}_{(\alpha,\beta)}^*(C,B)$ is a differential graded
  $k$--algebra, we must prove that the Leibniz rule holds.  Since
  $\C{K}_{(\alpha,\beta)}^*(C,B)$ as an algebra is generated by degree
  1 terms, it is enough to check the Leibniz rule for degree 0 and 1
  terms.  We see that
  \begin{align*}
    d(b'b)
     = & -(\B{I}\otimes b'b)
         + (\alpha(b'_{(1)})\alpha(b_{(1)})\B{I}\beta(b_{(3)})\beta(b'_{(3)})\otimes
            b'_{(2)}b_{(2)})\\
     = & -(\B{I}\otimes b'b) + (\alpha(b'_{(1)})\B{I}\beta(b'_{(3)})\otimes b'_{(2)}b)\\
       & - (\alpha(b'_{(1)})\B{I}\beta(b'_{(3)})\otimes b'_{(2)}b)
         + (\alpha(b'_{(1)})\alpha(b_{(1)})\B{I}\beta(b_{(3)})\beta(b'_{(3)})\otimes
            b'_{(2)}b_{(2)})\\
     = & d(b') b + b' d(b)
  \end{align*}
  for any $b,b'\in B$.  Moreover,
  \begin{align*}
    d((c\otimes b')b)
    = & d(c\otimes b'b)
    = -(\B{I}\otimes c\otimes b'b) + (c_{(1)}\otimes c_{(2)}\otimes b'b)
      -(c\otimes \alpha(b'_{(1)})\alpha(b_{(1)})\B{I}\beta(b_{(3)})\beta(b'_{(3)})
        \otimes b'_{(2)} b_{(2)})\\
    = & -(\B{I}\otimes c\otimes b'b) + (c_{(1)}\otimes c_{(2)}\otimes b'b)
      -(c\otimes \alpha(b'_{(1)}) \B{I} \beta(b'_{(3)}) \otimes b'_{(2)} b)\\
      & +(c\otimes \alpha(b'_{(1)}) \B{I} \beta(b'_{(3)}) \otimes b'_{(2)} b)
      -(c\otimes \alpha(b'_{(1)})\alpha(b_{(1)}) \B{I} \beta(b_{(3)}) \beta(b'_{(3)})
        \otimes b'_{(2)} b_{(2)})\\
    = & d(c\otimes b')b - (c\otimes b')d(b)
  \end{align*}
  for any $(c\otimes b')$ in $\C{K}_{(\alpha,\beta)}^1(C,B)$ and
  $b\in B$.  Then, again for the same elements
  \begin{align*}
   d(b(c\otimes b'))
    = & d(\alpha(b_{(1)}) c \beta(b_{(3)}) \otimes b_{(2)} b')\\
%    = &  -(\B{I}\otimes \alpha(b_{(1)}) c \beta(b_{(3)})\otimes b_{(2)} b')
%        + (\alpha(b_{(1)}) c_{(1)} \beta(b_{(5)}) \otimes
%           \alpha(b_{(2)}) c_{(2)} \beta(b_{(4)}) \otimes b_{(3)} b')\\
%      & -(\alpha(b_{(1)}) c \beta(b_{(3)}) \otimes
%          \alpha(b_{(2)(1)}) \alpha(b'_{(1)}) \B{I} \beta(b'_{(3)}) \beta(b_{(2)(3)})\otimes
%          b_{(2)(2)} b'_{(2)})\\
%      & -(\alpha(b_{(1)}) c \beta(b_{(5)}) \otimes
%          \alpha(b_{(2)})\alpha(b'_{(1)}) \B{I} \beta(b'_{(3)})\beta(b_{(4)})\otimes
%          b_{(3)} b'_{(2)})\\
    = &  -(\B{I}\otimes \alpha(b_{(1)}) c \beta(b_{(3)})\otimes b_{(2)} b')
         +(\alpha(b_{(1)})\B{I}\beta(b_{(5)})\otimes \alpha(b_{(2)}) c \beta(b_{(4)})\otimes
           b_{(3)} b')\\
      &  -(\alpha(b_{(1)})\B{I}\beta(b_{(5)})\otimes \alpha(b_{(2)}) c \beta(b_{(4)})\otimes
           b_{(3)} b')
        + (\alpha(b_{(1)}) c_{(1)} \beta(b_{(5)}) \otimes
           \alpha(b_{(2)}) c_{(2)} \beta(b_{(4)}) \otimes b_{(3)} b')\\
      & -(\alpha(b_{(1)}) c \beta(b_{(5)}) \otimes
          \alpha(b_{(2)})\alpha(b'_{(1)}) \B{I} \beta(b'_{(3)})\beta(b_{(4)})\otimes
          b_{(3)} b'_{(2)})\\
    = & d(b)(c\otimes b') + b d(c\otimes b')
  \end{align*}
  And finally for $(c\otimes 1)$ and $(c'\otimes 1)$ in
  $\C{K}_{(\alpha,\beta)}^*(C,B)$ we see that
  \begin{align*}
    d((c\otimes 1)(c'\otimes 1))
    = & d(c\otimes c'\otimes 1)\\
    = & -(\B{I}\otimes c\otimes c'\otimes 1)
       +(c_{(1)}\otimes c_{(2)}\otimes c'\otimes 1)
       -(c\otimes c'_{(1)}\otimes c'_{(2)}\otimes 1)
       +(c\otimes c'\otimes \B{I}\otimes 1)\\
    = & -(\B{I}\otimes c\otimes c'\otimes 1)
       +(c_{(1)}\otimes c_{(2)}\otimes c'\otimes 1)
       -(c\otimes\B{I}\otimes c'\otimes 1)\\
      &+(c\otimes\B{I}\otimes c'\otimes 1)
       -(c\otimes c'_{(1)}\otimes c'_{(2)}\otimes 1)
       +(c\otimes c'\otimes \B{I}\otimes 1)\\
    = & d(c\otimes 1)(c'\otimes 1) - (c\otimes 1) d(c'\otimes 1)
  \end{align*}
  Now, one can inductively prove that
  \begin{align*}
    d(\Psi\Phi) =  d(\Psi)\Phi + (-1)^{|\Psi|}\Psi d(\Phi)
  \end{align*}
  for any $\Psi$ and $\Phi$ in $\C{K}_{(\alpha,\beta)}^*(C,B)$ proving
  $\C{K}_{(\alpha,\beta)}^*(C,B)$ is a differential graded $k$--algebra.
\end{proof}

\begin{cor}
  $H_*(\C{K}_{(\alpha,\beta)}^*(C,B))$ is a graded algebra.
\end{cor}

Now we can  identify the homology of the
$(\alpha,\beta)$--equivariant differential calculus:
\begin{prop}
  $B$ is a $C$--comodule and
  $H_*(\C{K}_{(\alpha,\beta)}^*(C,B))$ is isomorphic to ${\rm
  Cotor}_C^*(k,B)$.
\end{prop}

\begin{proof}
  The $C$--comodule structure is given by
  $\rho_B(b)=\alpha(b_{(1)})\B{I}\beta(b_{(3)})\otimes b_{(2)}$ for
  any $b\in B$.  Now, one should observe that
  $\C{K}_{(\alpha,\beta)}^*(C,B)$ is the two sided cobar complex
  $\C{B}_*(k,C,B)$ of $C$ with coefficients in $k$ and $B$.
\end{proof}

\begin{defn}
  A $k$--module $X$ is called an $(\alpha,\beta)$--equivariant
  $C$--comodule if (i) $X$ is a left $C$--comodule via a structure
  morphism $\rho_X:X\to C\otimes X$ (ii) $X$ is a left $B$--module
   via a structure morphism $\mu_X:B\otimes X \to X$ (iii) one
  has
  \begin{align*}
    \rho_X(bx) = \alpha(b_{(1)})x_{(-1)}\beta(b_{(3)})\otimes b_{(2)}x_{(0)}
  \end{align*}
  for any $b\in B$ and $x\in X$.
\end{defn}

\begin{thm}
  The category of $(\alpha,\beta)$--equivariant $C$--comodules is
  equivalent to the category of $B$--modules admitting a flat
  connection with respect to the differential calculus
  $\C{K}_{(\alpha,\beta)}^*(C,B)$.
\end{thm}

\begin{proof}
  Assume we have a $k$--module morphism $\rho_X:X\to C\otimes X$ and
  define $\rho_X(x):=\nabla(x)+(1\otimes x)$ for any $x\in X$ where we
  denote $\rho_X(x)$ by $x_{(-1)}\otimes x_{(0)}$.  Now for any $b\in
  B$ and $x\in X$ we have
  \begin{align*}
    \nabla(bx)
     = & -(\B{I}\otimes bx) + (bx)_{(-1)}\otimes (bx)_{(0)}\\
     = & -(\B{I}\otimes bx) + (\alpha(b_{(1)})\B{I}\beta(b_{(3)})\otimes b_{(2)}x)
         - (\alpha(b_{(1)})\B{I}\beta(b_{(3)})\otimes b_{(2)}x)
         + (\alpha(b_{(1)})x_{(-1)}\beta(b_{(3)})\otimes b_{(2)}x_{(0)})\\
     = &\ d(b)\eotimes{B}x + b\nabla(x)
  \end{align*}
  iff $(bx)_{(-1)}\otimes (bx)_{(0)} =
  \alpha(b_{(1)})x_{(-1)}\beta(b_{(3)})\otimes b_{(2)}x_{(0)}$.  In
  other words $\nabla:X\to C\otimes X$ is a connection iff the
  morphism of $k$--modules $\rho_X:X\to C\otimes X$ satisfies the 
  $(\alpha,\beta)$--equivariance condition.  Moreover, if we extend
  $\nabla$ to $\widehat{\nabla}:C\otimes X\to C\otimes C\otimes X$
  by letting
  \begin{align*}
    \widehat{\nabla}(c\otimes x)
      = & d(c\otimes 1)\eotimes{B}x - (c\otimes 1)\eotimes{B}d(x)
  \end{align*}
  for any $(c\otimes x)\in C\otimes X$, then  we have
  \begin{align*}
    \widehat{\nabla}^2(x)
    = & - d(\B{I}\otimes 1)\eotimes{B}x + (\B{I}\otimes 1)\eotimes{B}\nabla(x)
        + d(x_{(-1)}\otimes 1)\eotimes{B} x_{(0)}
        - (x_{(-1)}\otimes 1)\eotimes{B}\nabla(x_{(0)})\\
    = & (\B{I}\otimes\B{I}\otimes x) - (\B{I}\otimes\B{I}\otimes x)
        + (\B{I}\otimes\B{I}\otimes x) - (\B{I}\otimes\B{I}\otimes x)
        +(\B{I}\otimes x_{(-1)}\otimes x_{(0)})\\
      & -(\B{I}\otimes x_{(-1)}\otimes x_{(0)})
        +(x_{(-1)(1)}\otimes x_{(-1)(2)}\otimes x_{(0)})
        -(x_{(-1)}\otimes \B{I}\otimes x_{(0)})\\
      & +(x_{(-1)}\otimes \B{I}\otimes x_{(0)})
        -(x_{(-1)}\otimes x_{(0)(-1)}\otimes x_{(0)(0)})\\
    = & (x_{(-1)(1)}\otimes x_{(-1)(2)}\otimes x_{(0)})
        -(x_{(-1)}\otimes x_{(0)(-1)}\otimes x_{(0)(0)}) = 0
  \end{align*}
  iff $\rho_X:X\to C\otimes X$ is a coassociative coaction of $C$.
  The result follows.
\end{proof}

\begin{defn}
  Let $X$ be an $(\alpha,\beta)$--equivariant $C$--comodule, i.e. $X$
  admits a flat connection $\nabla:X\to C\otimes X$ with respect to
  the differential calculus $\C{K}_{(\alpha,\beta)}^*(C,B)$.  Define
  $\C{K}_{(\alpha,\beta)}^*(C,B,X)$ as
  $\C{K}_{(\alpha,\beta)}^*(C,B)\eotimes{B}X$ with the extended
  connection  $\widehat{\nabla}$ as its differential.
\end{defn}

\begin{prop}
  For any $(\alpha,\beta)$--equivariant module $X$  one has $H_*(\C{K}_{(\alpha,\beta)}^*(C,B,X))\cong {\rm
  Cotor}_C^*(k,X)$ where we think of $k$ as a $C$--comodule via the
  grouplike element $\B{I}:k\to C$.
\end{prop}

\begin{rem}
  Assume $B$ is a Hopf algebra. If $\alpha=id_B$,  $\beta=S$, and
  $C=B$, then  the differential calculus is $\widehat{\mathcal{K}}(H)$  and
  $(\alpha,\beta)$--equivariant $C$--comodules are Yetter-Drinfeld
  modules.  In case $\alpha=id_B$,  $\beta=S^{-1}$, and $C=B$, then  the
  differential calculus is ${\mathcal{K}}(H)$  and
  $(\alpha,\beta)$--equivariant $C$--comodules are
  anti-Yetter-Drinfeld modules.
\end{rem}

% \bibliographystyle{plain}
% \bibliography{bibliography}

\end{document}